\documentclass{singlecol-new}

\usepackage{natbib,stfloats}
\usepackage{mathrsfs}
\usepackage{amssymb}
\usepackage{amsmath}
\usepackage{epsfig}
\usepackage{multirow}
\usepackage{longtable}
\usepackage{algorithmic}
\usepackage{lscape}

\def\0{\bf \0}
\def\A{{\bf A}}

\def\I{{\bf I}}

\def\L{{\bf L}}
\def\M{{\bf M}}

\def\0{{\bf 0}}

\def\R{\mathbb{R}}

\def\T{{\bf T}}
\def\U{{\bf U}}

\def\a{{\bf a}}
\def\b{{\bf b}}
\def\c{{\bf c}}

\def\p{{\bf p}}
\def\q{{\bf q}}

\def\x{{\bf x}}
\def\y{{\bf y}}

\def\Tr{{\rm T}}
\def\T{{\rm T}}

\theoremstyle{TH}{

\newtheorem{algorithm}{Algorithm}
}

\theoremstyle{THrm}{

\newtheorem{remark}{Remark}

}

\theoremstyle{THhit}{

}

\makeatletter
\def\theequation{\arabic{equation}}

\makeatother

\begin{document}%

\setcounter{page}{1}

\LRH{Yang and Vitor}

\RRH{A Double-Pivot Degenerate-Robust Simplex Algorithm for Linear Programming}

\VOL{x}

\ISSUE{x}

\PUBYEAR{202X}

\BottomCatch

\CLline

\subtitle{}

\title{A Double-Pivot Degenerate-Robust Simplex Algorithm for Linear Programming}

\authorA{Yaguang Yang}

\affA{U.S. Department of Commerce,\\ Rockville, MD, USA \\
E-mail: yaguang.yang@verizon.net}

\authorB{Fabio Vitor}
\affB{Department of Mathematics,\\ University of Nebraska at Omaha,\\ Omaha, NE, USA \\
E-mail: ftorresvitor@unomaha.edu}

\begin{abstract}
A double pivot simplex algorithm that combines features of two recently published papers by these authors is proposed. The proposed algorithm is implemented in MATLAB. The MATLAB implementation is tested, along with a MATLAB implemention of Dantzig's algorithm, for several test sets, including a set of cycling linear programming problems, Klee-Minty's problems, randomly generated linear programs, and Netlib benchmark problems. The test results show that the proposed algorithm, with a careful implementation, is (a) degenerate-robust as expected, and (b) more efficient than Dantzig's algorithm for large size randomly generated linear programming problems, but less efficient for Netlib benchmark problems and small size randomly generated problems in terms of CPU time.
\end{abstract}

\KEYWORD{Double Pivots; Degenerate-Robust; Simplex Method; Linear Programming; Klee-Minty Cube.}

\REF{to this paper should be made as follows: Yang, Y. and Vitor, F. (202X) `A double-pivot degenerate-robust simplex algorithm for linear programming', {\it International Journal of Operational Research}, Vol. x, No. x, pp.xxx\textendash xxx.}

\begin{bio}
Yaguang Yang received his B.S. and M.S. degrees from Huazhong University of Science and Technology, China. From 1985 to 1990, he was a Lecturer at Zhejiang University in China. In 1996, he received his PhD degree from the Department of Electrical and Computer Engineering at the University of Maryland, College Park. He has been working on linear optimization, nonlinear optimization, optimization on Riemannian manifolds, interior-point methods, and their applications in engineering problems. Currently, he works for the U.S. Department of Commerce.\vs{9}

\noindent Fabio Vitor is an Assistant Professor in the Department of Mathematics at
the University of Nebraska at Omaha. He received his PhD in Industrial
Engineering and MS in Operations Research from Kansas State University,
and his BS in Industrial Engineering from Maua Institute of Technology,
Brazil. He also worked for Monsanto, Kalmar (Cargotec Corporation), and
Volkswagen. His research interests include applied optimization topics and the development of algorithms to more quickly solve continuous and discrete optimization problems such as linear, nonlinear, and integer programs.\vs{8}
\end{bio}

\maketitle

\section{Introduction} \label{Intro}

It has been more than 70 years since Dantzig formulated
a linear programming (LP) problem and proposed the simplex 
method \citep{dantzig49}. The main idea of this algorithm
is to search for an optimizer from a vertex to the next vertex in
the polyhedron formed by the linear constraints using
Dantzig's pivoting rule. Since then, linear programming 
has been one of the mostly studied problems in operations research \citep{borgwardt20, bowly20, dadush20, cavaleiro21, suriyanarayana22}.
The simplex method remained to be a mainstream technique
for LP until \cite{km72} found an example showing that, in the worst case scenario, Dantzig's pivoting rule
needs exponentially many iterations to find an optimal
solution of the LP problem. Klee and Minty's work inspired 
researchers to find different ways to solve LP problems.

Interior-point methods emerged as a computationally
feasible alternative technique, which admits algorithms
that find an optimal solution in finite iterations bounded
by a polynomial of the problem size. In the last three 
decades, most researches in linear programming 
focused on finding novel interior-point methods
\citep{wright97,gondzio12,santos19,yang20a,ve22}. However, the simplex method
never faded away. Motivated by enhancing
the computational performance of the simplex method, \cite{ve18} proposed a double pivoting algorithm, 
which updates two pivot variables at a time using 
Dantzig's pivoting rule. Motivated by improving the bound of iteration number
of the simplex method, \cite{yang20} independently
proposed a different double pivoting rule that uses
a combined criteria to select two pivot variables.

Realizing the merits of the pivoting rule proposed by \cite{yang20} and 
the slope algorithm created by \cite{ve18}, this paper proposes and implements a new double pivoting algorithm that combines the pivoting rule of Yang and
the slope algorithm designed to find an optimal solution and an optimal basis of two-dimensional linear programming problems. The proposed algorithm selects two entering variables based on three criteria. The first entering variable is selected based on Dantzig's pivoting rule, the second entering variable is selected based on the longest step size rule, the combination
of the coefficients of these two variables are determined by the criterion that optimizes the cost reduction, which is equivalent to solving a linear programming
problem of two variables. Notice that the slope algorithm was developed specifically to solve this type of problem.

The intuitions behind the pivoting criteria of Yang are
as follows: first, Dantzig's pivoting rule is selected because
the computational experience in decades shows that 
it is a cheap pivoting rule compared to other
popular pivoting rules \citep{ps14}; second, specially 
designed examples show that none of the popular
pivoting rules is better than others \citep{psz08} in the
worst case, hence, using a combination of different rules in a 
random way can be beneficial and has been proved that the
strategy results in some polynomial simplex algorithms in the 
sense of statistics \citep{ghz98, ks06};  third, using the longest
step size rule will decrease the chance to enter a degenerate
basic feasible solution and increase the chances to get
out of a degenerate solution \citep{yang21}; fourth, the 
longest step rule may be beneficial to reduce the 
iteration numbers \citep{yang20}; finally, using deterministic
pivoting rules may give some hope to find a strongly
polynomial algorithm to solve LP problems, which is not possible for 
randomized methods such as the ones proposed by \cite{ghz98} and \cite{ks06}.

Observe that another popular pivoting rule is Harris steepest-edge rule \citep{fg92}, which is similar to the longest step rule in the sense that it needs to compute the step size for all possible pivoting candidates corresponding to a negative reduced cost. Several computational experiments showed that the steepest-edge rule needs fewer simplex iterations than the most negative reduced cost rule \citep{gr77, ps14} as one would expect. However, computational experience shows	that finding all step sizes for possible pivoting candidates corresponding to the negative coefficients in the revised cost function can be extremely expensive if the cardinality of the negative reduced costs is large. For more complete surveys on pivoting rules, see \cite{terlaky93}, \cite{pan14}, and \cite{ploskas17}.

Throughout the paper, capitalized bold font is used for matrices, small case bold font is used for vectors, and normal font is used for scalars. To save space, stacked vectors $\left[ \x^{\T}, \y^{\T} \right]^{\T}$ are represented as $(\x, \y)$. The remainder of the paper is organized as follows. Section \ref{Prop_Algo} presents the proposed algorithm. Section \ref{Imp_Det} discusses some implementation details. Numerical test problems and test results are provided in Section \ref{Num_Tests}. Conclusion remarks are summarized in Section \ref{Conc}.

\section{Proposed Algorithm} \label{Prop_Algo}

Consider the primal linear programming problem in standard form:

\begin{eqnarray}
\begin{array}{ll}
\min & \ \ \c^{\Tr}\x, \\ 
\mbox{\rm subject to} & \ \ \A\x=\b, \hspace{0.1in} \x \ge \0,
\end{array}
\label{LP}
\end{eqnarray}

\noindent where $\A \in {\R}^{m \times n}$, $\b \in {\R}^{m}$, $\c \in {\R}^{n}$ are given, and $\x \in {\R}^n$  is the vector to be optimized. For the following discussion, two primary assumptions are made:

\begin{enumerate}
\item $\mbox{rank}(\A) = m$; \label{assump_rank}
\item primal LP problem (\ref{LP}) has an optimal solution. \label{assump_opt}
\end{enumerate}

These assumptions are standard. A feasible solution of the linear program satisfies
the conditions of $\A \x =\b$ and $\x \ge \0$, which exists because of Assumption \ref{assump_opt}. Among all feasible solutions, this paper considers only basic feasible solutions, which corresponds to the vertices of the convex polytope
described by the constraints of (\ref{LP}). Denote by $B \subset \{ 1, 2, \ldots, n \}$ the index set with cardinality $|B| = m$ and $N = \{ 1, 2, \ldots, n \} \setminus B$ the complementary set of $B$ with cardinality $|N|=n-m$ such that 
matrix $\A$ and vector $\x$ can be partitioned as $\A=[ \A_B, \A_N]$ and $\x =(\x_B, \x_N)$ with $\x_B \ge 0$ and $\x_N=0$. Moreover, if $\A_B \x_B = \b$, then $\x=(\x_B,\0)$ is a basic feasible solution. Because some components of $\x_B$ may be zeros, and/or even worse, $\A_B$ may not be full rank, this partition admits a degenerate basic feasible solution, which are seen in many real-world applications and also in many Netlib benchmark problems \citep{bdgr95}. Denote by $\mathcal{B}$ the set of all bases $B$ and by $\mathcal{N}$ the set of all non-bases $N$. Therefore, the linear programming problem (\ref{LP}) can be written as:

\begin{eqnarray}
\begin{array}{ll}
\min & \ \ \c_B^{\T}\x_B + \c_N^{\T}\x_N,  \\
\mbox{\rm subject to} 
& \ \ \A_B\x_B + \A_N\x_N=\b, 
\hspace{0.1in} \x_B \ge \0, \hspace{0.1in}  \x_N \ge \0.
\end{array}
\label{LP1}
\end{eqnarray}

Let superscript $k$ represent the $k$th iteration. Thus, the matrices and vectors in the $k$th iteration are denoted by $\A_{B^k}$, $\A_{N^k}$, $\c_{B^k}$, $\c_{N^k}$, $\x_{B^k}$, and $\x_{N^k}$ where $\x^k=(\x_{B^k},\x_{N^k})$ is a basic feasible solution of (\ref{LP}) with $\x_{B^k} \ge \0$ and $\x_{N^k}= \0$. Similarly, define $\x^*= (\x_{B^*}, \x_{N^*})$ as an optimal basic feasible solution of (\ref{LP}) with $\A_{B^*}\x_{B^*}=\b$, $\x_{B^*} \ge \0$, $\x_{N^*}=\0$, and $z^* = \c^{\Tr} \x^*$ as the optimal objective function value. Observe that the 
partition of $(B^k,N^k)$ keeps updating and it is different from the partition $(B^*,N^*)$ before an optimizer is found. If $\A_{B^k}$ is full rank, the
reduced cost vector can be calculated as

\begin{equation}
\bar{\c}_{N^k}^{\Tr} =( \c_{N^k} - \A_{N^k}^{\Tr} 
\A_{B^k}^{-\Tr} \c_{B^k} )^{\Tr},
\label{cNbar}
\end{equation}

\noindent and the iterate of $\x^k=\A_{B^k}^{-1}\b$.

\subsection{LU Decomposition versus Pseudo Inverse Solution}
\label{avoid}

Avoiding the calculation of the inverse is not only useful to reduce the computational cost, but also imperative to make sure that the algorithm works when 
$\A_B$ is singular. Notice that the selection of the entering variable in the simplex method does not consider if the new $\A_B$ is full rank or not. As a matter of fact, the authors noticed that for Netlib benchmark problems, it is not unusual that $\A_B$ is not full rank. For this paper's implementation, LU decomposition is used to achieve this. Let

\begin{equation}
\A_{B^k} =\L \U,
\label{LUdecomp}
\end{equation}

\noindent where $\L$ is a full rank permuted lower triangular matrix and $\U$ is an upper triangular matrix with the same rank of $\A_{B^k}$. Therefore, if 
$\A_{B^k}$ is not full rank, $\U$ is not full rank, i.e., some diagonal elements of $\U$ are zeros. Given the LU decomposition (\ref{LUdecomp}), the reduced cost vector can be calculated as follows. Let $\bar{\p}=\A_{B^k}^{-\Tr} \c_{B^k}$, i.e.,
$\A_{B^k}^{\Tr}\bar{\p}= \c_{B^k}$. Thus, $\bar{\p}$ is obtained by sequentially
solving two systems of linear equations

\begin{equation}
\U^{\Tr} \p=\c_{B^k} \hspace{0.1in} \text{and}
\hspace{0.1in} \L^{\Tr} \bar{\p}=\p.
\label{LU1}
\end{equation}

\noindent Hence, the reduced cost vector is obtained by computing

\begin{equation}
\bar{\c}_{N^k}=\c_{N^k}- \A_{N^k}^{\Tr} \bar{\p}.
\label{reducedCost}
\end{equation}

\noindent The iterate $\x^k=\A_{B^k}^{-1}\b$ can be obtained by sequentially solving two systems of linear equations as follows:

\begin{equation}
\L \q=\b \hspace{0.1in} \text{and} \hspace{0.1in} \U \x^k=\q.
\label{LU2}
\end{equation}

When $\A_{B^k}$ is not full rank, and some diagonal elements of $\U$ are zeros, one may consider the use of the pseudo inverse of $\A_{B^k}$. A simple and careful 
analysis shows that this would not work in some cases. The pseudo inverse is equivalent to the following process. Let $\bar{\U}$ be the matrix that removes
rows with zero diagonal elements from $\U$ and $\bar{\q}$ be the vector that removes the corresponding elements in $\q$; then Equation (\ref{LU2}) becomes $\bar{\U}\x^k =\bar{\q}$ and the solution is given by

\begin{equation}
\x^k = \begin{cases}
\bar{\U}^{\Tr}(\bar{\U}\bar{\U}^{\Tr})^{-1}\bar{\q}
& \ \ {\mbox {if}} ~~ \bar{\q} \neq 0, \\
\x^k \in  \bar{\U}^{\perp}
&  \ \ {\mbox{otherwise.}}  \end{cases}
\label{LU2b}
\end{equation}

The matrix $\bar{\U}^{\perp}$ can be obtained by using QR decomposition for $\bar{\U}$. One can easily see that the least squared solution \citep{bg03} of (\ref{LU2b}) may not meet the equations corresponding to the deleted rows of $\U$.
Testing on some Netlib problems verified the analysis. The strategy used in this paper is to add a small $\epsilon$ to the zero diagonal elements of $\U$ when it is singular. This will avoid the difficulty in solving (\ref{LU1}) and (\ref{LU2}).

\subsection{A Double Pivot Algorithm}

The double pivot algorithm in this paper is based on \cite{yang20}. First, how the entering variables are selected is briefly describe. Let $C^k$ be the index set defined as follows:

\begin{equation}
C^k \in \{ j^k ~|~ \bar{c}_{j^k} <0 \},
\label{RcostIndex}
\end{equation}

\noindent and the cardinality $ | C^k | =p$. Clearly, if $\bar{\c}_{N^k} \ge \0$, then $C^k = \emptyset$, and an optimizer is found. If $\bar{c}_{j^k} <0$ for some 
$j^k \in C^k$, then an entering variable $x_{\jmath^k}$ in the next vertex is chosen from the set of $\{ j^k ~|~ \bar{c}_{j^k} <0 \}$ because by increasing $x_{\jmath^k}$, the objective function $\c^{\Tr} \x=\c_{B^k}^{\T}\x_{B^k}+\bar{c}_{\jmath^k}x_{\jmath^k}$ will be reduced. The first entering variable $x_{\jmath_1^k}$ is selected by using Dantzig's pivoting rule:
 
\begin{equation}
{\jmath_1^k} :=\{ {\jmath_1^k} ~|~ \bar{c}_{\jmath_1^k} 
= \min_{j^k \in C^k} \bar{\c}_{j^k}  \}.
\label{enteringVar}
\end{equation} 

\noindent If there is a tie, the minimum index of ${\jmath_1^k}$ will be used to break the tie. Denote $\bar{\b} = \A_{B^k}^{-1} \b$ and $\bar{\a}_{\jmath^k} = \A_{B^k}^{-1} \A_{\jmath^k}$, where $\jmath^k$ is in $C^k$ and determined by (\ref{enteringVar}). Also, denote by $\bar{b}_i$ the $i$th element of $\bar{b}$ and by $\bar{a}_{j^k,i}$ the $i$th element of $\bar{\a}_{\jmath^k}$. Thus, the leaving variable for Dantzig's pivoting rule is to select $x_{\imath_1^k}$ that satisfies

\begin{equation}
x_{\imath_1^k} = 
\min_{i \in  \{ 1,\ldots,m \}} \bar{b}_i / \bar{a}_{{\jmath_1^k},i},
\hspace{0.1in} \mbox{subject to}
\hspace{0.1in}    \bar{a}_{{\jmath_1^k},i}>0.
\label{iniVar2}
\end{equation}

\noindent If there is a tie, the minimum index of ${\imath_1^k}$ will
be used to break the tie. The reduced cost value is given by $\bar{c}_{\imath_1^k} x_{\imath_1^k}$.

If the cardinality $|C^k| \ge 2$, the second entering variable 
whose index is $\jmath_2^k$ will maximize the step-size, i.e.,

\begin{equation}
x_{\jmath_2^k} = \max_{\bar{c}_{j^k}<0,
\bar{c}_{j^k}\neq\bar{c}_{\jmath_1^k}}
\Big\{
\min_{i \in  \{ 1,\ldots,m \}} \bar{b}_i / \bar{a}_{j^k,i},
\hspace{0.1in} \mbox{subject to}
\hspace{0.1in}    \bar{a}_{j^k,i}>0
\Big\}.
\label{iniVar3}
\end{equation}

The leaving variable $x_{\imath_2^k}$ is determined by the index that achieves the max-min value in (\ref{iniVar3}). Given the entering and leaving variables determined by the longest step size rule, one can determine the next iterate $\x^{k+1}$ and the reduced cost $\bar{c}_{\imath_2^k} x_{\imath_2^k}$.

This paper discusses two special cases of the double pivot algorithm. For the general case, given the two entering variables $x_{\jmath_1^k}$ and $x_{\jmath_2^k}$, one needs to determine the two leaving variables which is equivalent to solving a linear programming problem with constraints in a two-dimensional space. Let $\A_{(\jmath_1,\jmath_2)}$ be composed of the 
$\jmath_1$ and $\jmath_2$ columns of $\A_N$, and $\jmath_1$ and $\jmath_2$ be determined by (\ref{enteringVar}) and (\ref{iniVar3}). Let $\bar{\A}_{(\jmath_1,\jmath_2)} = \A_B^{-1} \A_{(\jmath_1,\jmath_2)}$ and $\bar{\c}_{(\jmath_1,\jmath_2)} < \0$ be the two corresponding elements in $\bar{\c}_N$. For the two entering indexes $(\jmath_1,\jmath_2) \in C^k$ satisfying
$\x_{(\jmath_1,\jmath_2)} =(x_{\jmath_1},x_{\jmath_2}) \ge \0$, one needs

\begin{equation}
\x_{B^k}^{k+1}=\A_{B^k}^{-1}\b
-\A_{B^k}^{-1}\A_{N^k}\x_{N^k} 
= \bar{\b} - \bar{\A}_{(\jmath_1,\jmath_2)}
\x_{(\jmath_1,\jmath_2)} \ge \0.
\label{updatedXB}
\end{equation}

Observe that instead of using matrix inverse, $ \bar{\b}$ and $\bar{\A}_{(\jmath_1,\jmath_2)}$ are obtained by using LU decomposition and then solving the systems of linear equations discussed in Section \ref{avoid}. Therefore, the problem of finding a new vertex is reduced to minimizing the following linear programming problem:

\begin{eqnarray}
\begin{array}{ll}
\min  & \ \ \bar{\c}_{(\jmath_1,\jmath_2)}^{\T}
\x_{(\jmath_1,\jmath_2)}, \\ 
\mbox{\rm subject to} 
& \ \ \bar{\A}_{(\jmath_1,\jmath_2)}
\x_{(\jmath_1,\jmath_2)} \le \bar{\b}, 
\hspace{0.1in} \x_{(\jmath_1,\jmath_2)} \ge \0.
\end{array}
\label{twoDim}
\end{eqnarray}

\noindent Again, the vector $\bar{\c}_{(\jmath_1,\jmath_2)}$ is obtained by using LU decomposition and solving the systems of linear equations discussed in Section \ref{avoid}. The third merit criterion is then introduced, which is to
determine the values of the two entering variables to minimize the objective function under the constraints of (\ref{twoDim}). Remarks \ref{polynomial} and \ref{cycling} formalize the pivoting strategy of the proposed algorithm. 

\begin{remark}
Using multiple merit criteria to select entering variables has been proved to be an effective strategy to solve LP problems and has been used in randomized pivoting algorithms \citep{ghz98, ks06}. Since almost all popular deterministic simplex pivoting rules are proved not polynomial \citep{friedmann11, gs79, jeroslow73, psz08}, using a combined deterministic merit criteria to select entering variables may give some fresh ideas to find a strongly polynomial algorithm (defined in \cite{smale99}) to solve linear programming problems. \label{polynomial}
\end{remark}

\begin{remark}
Using the longest step size rule to select an entering parameter will increase the chance to have a non-degenerate basic feasible solution because the entering variable is more likely to be greater than zero than any other criteria. Since the double pivoting rule will optimize the cost function, the step-sizes for the two pivot variables will not be all zeros if the longest step size is not zero. Therefore, the authors claim that the double pivot algorithm is degenerate-robust. This claim has been observed in numerical experiments to be discussed later. \label{cycling}
\end{remark}

Even though problem (\ref{twoDim}) can be solved using standard simplex algorithms or interior-point methods, there is a much more efficient technique \citep{ve18, vitor18, vitor19, vitor22} for this special LP problem, which is discussed in the next section.

\subsection{The Slope Algorithm to Solve Two-Variable Linear Programs}

The slope algorithm is a fast technique that runs in O$(m \log m)$ time and can find both an optimal solution and an optimal basis of two-variable linear programs. The slope algorithm was first proposed by \cite{ve18} and later improved by \cite{vitor18, vitor19, vitor22}. Observe that the slope algorithm can quickly solve problem (\ref{twoDim}). Formally, define a two-variable linear program (2VLP) as:

\begin{eqnarray}
\begin{array}{ll}
\max  &  \c'^{\Tr}\x', \\ 
\mbox{\rm subject to} 
&  \A'\x'\leq\b',~~~~\x' \ge \0,
\end{array}
\label{2VLP}
\end{eqnarray}

\noindent where $\A' \in {\R}^{(m+2) \times 2}$, $\b' \in {\R}^{m+2}$, 
$\c' \in {\R}^{2}$ are given, and $\x' \in {\R}^2$  is the 
vector of decision variables. For a 2VLP, $A'_{m+1,1} = A'_{m+2,2} = -1$, 
$A'_{m+1,2} = A'_{m+2,1} = 0$, and $b'_{m+1} = b'_{m+2} = 0$. 
Observe that these two constraints are exactly the nonnegative conditions of 
(\ref{twoDim}). For simplicity, the slope algorithm assumes that $c'_1 > 0$, 
$c'_2 > 0$, and $\b'_i \geq 0$ for each $i \in \{1,2,...,m,m+1,m+2\}$. 
One can easily see that problem (\ref{twoDim}) satisfies these assumptions.

Overall, the slope algorithm contrasts the ``slope'' formed by the objective 
function coefficients of the 2VLP, $c'_1$ and $c'_2$, with the ``slope'' of 
every constraint in the problem. The slope of each constraint 
$i \in \{1,2,...,m,m+1,m+2\}$ is defined by $\alpha_i$ and can be computed as:

\begin{equation}
\alpha_i = 
\begin{cases}
-2M & \text{if $A'_{i,1}=0 \ \ \ \mbox{and} \ \ \ A'_{i,2}<0$}\\
-M + \frac{A'_{i,2}}{A'_{i,1}} & \text{if $A'_{i,1}>0 \ \ \ \mbox{and} \ \ \ A'_{i,2}<0$}\\
-M & \text{if $A'_{i,1} > 0 \ \ \ \mbox{and} \ \ \ A'_{i,2} = 0$}\\
\frac{A'_{i,2}}{A'_{i,1}} & \text{if $A'_{i,1} > 0 \ \ \ \mbox{and} \ \ \ A'_{i,2} > 0$}\\
M & \text{if $A'_{i,1} = 0 \ \ \ \mbox{and} \ \ \ A'_{i,2} > 0$}\\
M - \frac{A'_{i,1}}{A'_{i,2}} & \text{if $A'_{i,1} < 0 \ \ \ \mbox{and} \ \ \ A'_{i,2} > 0$}\\
2M & \text{if $A'_{i,1} < 0 \ \ \ \mbox{and} \ \ \ A'_{i,2} = 0$}\\
3M & \text{if $A'_{i,1} = 0 \ \ \ \mbox{and} \ \ \ A'_{i,2} = 0$}\\
3M & \text{if $A'_{i,1} < 0 \ \ \ \mbox{and} \ \ \ A'_{i,2} < 0$}\\
\end{cases} \label{eq_alpha}
\end{equation}

\noindent where

\begin{equation}
M > \max{\Bigg\{M', M'', \frac{c'_2}{c'_1}\Bigg\}} \label{eq_M}
\end{equation}

\noindent is a sufficiently large positive number and

\begin{equation}
M' = \max\Bigg\{\Bigg|\frac{A'_{i,1}}{A'_{i,2}}\Bigg| : 
A'_{i,2} \neq 0 \ \ \forall i \in \{1,2,...,m,m+1,m+2\}\Bigg\}
\end{equation}

\begin{equation}
M'' = \max\Bigg\{\Bigg|\frac{A'_{i,2}}{A'_{i,1}}\Bigg| : 
A'_{i,1} \neq 0 \ \ \forall i \in \{1,2,...,m,m+1,m+2\}\Bigg\}.
\end{equation}

\begin{figure*}[h]
\caption{Eight out of Nine Types of Constraints in a 2VLP and their Corresponding $\alpha$ Values} \label{Slopes}
\centerline{\epsffile{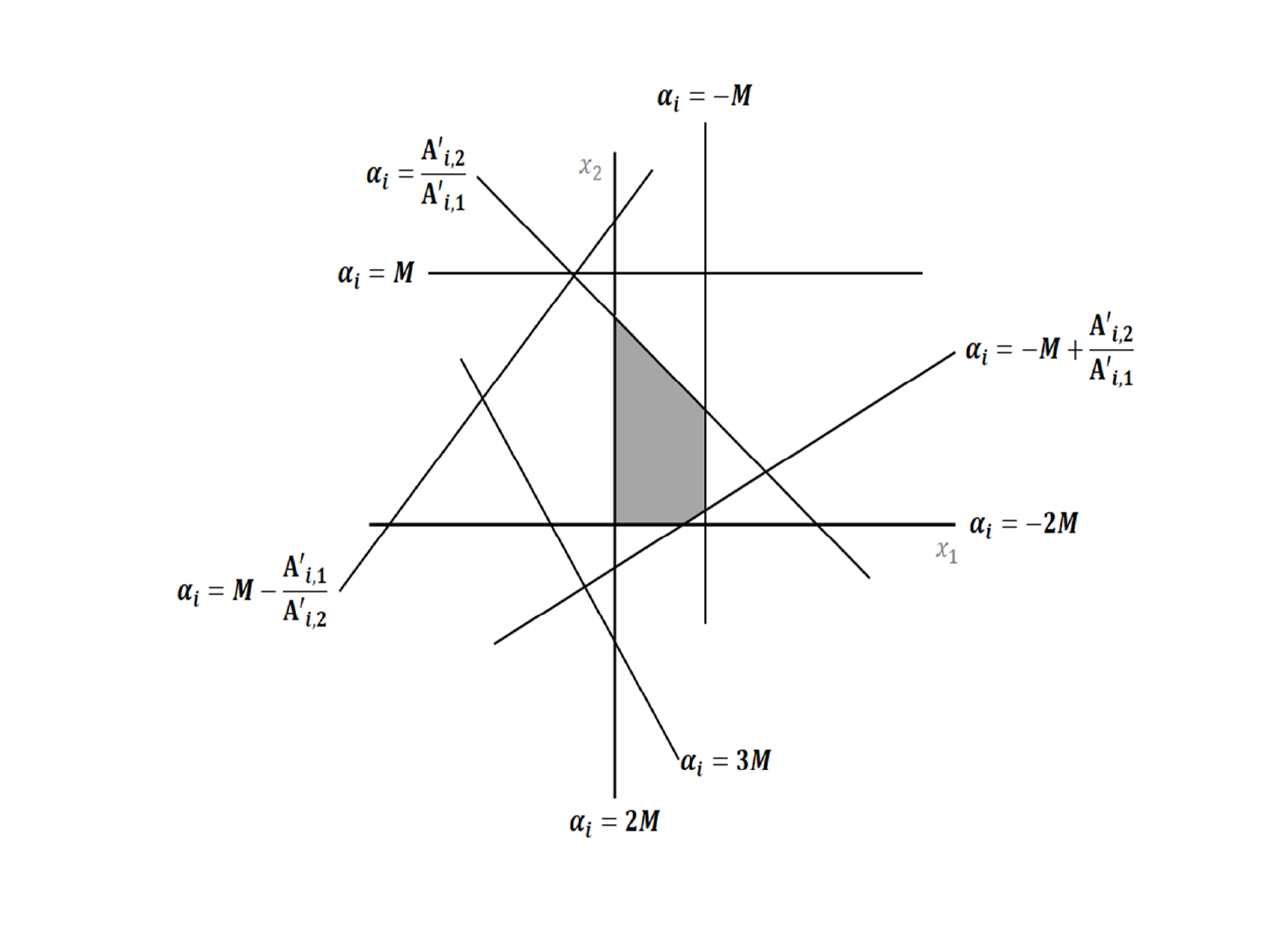}}
\end{figure*}

Notice that $M$ is used to determine some not well-defined slopes 
and differentiate the order of constraints. Figure \ref{Slopes} shows 
eight out of the nine types of constraints in a 2VLP (except 
constraints where $A'_{i,1} = 0$ and $A'_{i,2} = 0$ since those 
define the entire two-dimensional space) and their corresponding 
$\alpha$ values. Furthermore, the reader may observe that viewing 
the constraints in ascending order of the $\alpha$ values defines a 
counterclockwise orientation of the constraints.

Algorithm \ref{SAlgorithm} depicts the slope algorithm step by step. 
The input to Algorithm \ref{SAlgorithm} is a 2VLP and the method begins 
by computing a sufficiently large positive number $M$, calculating the 
``slope'' $\alpha_i$ of each constraint, and sorting the constraints in 
ascending order according to $\alpha_i$. From there, the slope algorithm 
finds two constraints, $\eta_{j'}$ and $\eta_{k'}$, which slope falls in 
between the slope of the cost coefficients $c'_1$ and $c'_2$ based on 
the sorted order. The method then checks whether the given 2VLP is 
unbounded or not. This can be done by the checking the presence/absence 
of some specific constraints that define a ray of unboundedness. That is, 
2VLP is unbounded if:

\begin{equation}
\begin{gathered}
\alpha_{\eta_{j'}} = -2M \ \ \text{and} \ \ \alpha_{\eta_{k'}} \geq M, \ \ \text{or}\\
-2M < \alpha_{\eta_{j'}} < -M \ \ \text{and} \ \ \alpha_{\eta_{k'}} = 2M, \ \ \text{or}\\
\alpha_{\eta_{j'}} = -M \ \ \text{and} \ \ \alpha_{\eta_{k'}} = 2M, \ \ \text{or}\\
-2M < \alpha_{\eta_{j'}} < -M \ \ \text{and} \ \ M < \alpha_{\eta_{k'}} < 2M, \ \ \text{and} \ \ \frac{A'_{\eta_{j'},2}}{A'_{\eta_{j'},1}} \leq \frac{A'_{\eta_{k'},2}}{A'_{\eta_{k'},1}}.\\
\label{eq_unb}
\end{gathered}
\end{equation}

If the given 2VLP is bounded, the slope algorithm continues and finds 
whether the intersection of constraints $\eta_{j'}$ and $\eta_{k'}$ is 
feasible on every other constraint. If not feasible, the algorithm replaces 
one of the two constraints with the constraint that violates the feasibility 
check, and repeat the process until a pair of feasible constraints is found. 
When this process is completed, the slope algorithm returns the optimal 
solution and also the two constraints that intersect at the optimal basis. 
That is, an optimal simplex basic feasible solution from where there does 
not exist a feasible improving search direction. This is possible because 
the slope algorithm finds a pair of constraints such that 
$\alpha_{\eta_{j'}} < \frac{c'_2}{c'_1} \leq \alpha_{\eta_{k'}}$ and 
$\alpha_{\eta_{k'}} - \alpha_{\eta_{j'}}$ is minimized. The reader is
encouraged to see \cite{ve18} and \cite{vitor18, vitor19, vitor22} for additional theoretical results.

\begin{algorithm} {The Slope Algorithm} \\
\begin{algorithmic}[1] 
\STATE Data: Matrix $\A'$, vectors $\b'$ and $\c'$
\STATE Compute a sufficiently large positive number $M$ according to (\ref{eq_M})
\STATE Compute $\alpha_i$ for each constraint 
$i \in \{1,2,...,m,m+1,m+2\}$ according to (\ref{eq_alpha})
\STATE Let $H = (\eta_1,\eta_2,...,\eta_m,\eta_{m+1},\eta_{m+2})$ 
be the list of constraint indices sorted in ascending order according to their $\alpha$ values
\STATE Find constraints $j'$ and $k' \in \{1,2,...,m,m+1,m+2\}$ 
such that $\alpha_{\eta_{j'}} < \frac{c'_2}{c'_1} \leq \alpha_{\eta_{k'}}$
\IF{(\ref{eq_unb}) is satisfied}
\STATE Report 2VLP is unbounded
\ELSE
\STATE $j \leftarrow j'$
\STATE $k \leftarrow k'$
\STATE Find the intersection of constraints $\eta_{j'}$ and
$\eta_{k'}$ and let its solution be $\bar{\x} = (\bar{x}_1,\bar{x}_2)$
\WHILE{$j > 1$ or $k < m+2$}
\IF{$j > 1$}
\STATE $j \leftarrow j-1$
\ENDIF
\IF{$A'_{\eta_j,1}\bar{x}_1 + A'_{\eta_j,2}\bar{x}_2 > b'_{\eta_j}$}
\STATE $j' \leftarrow j$
\STATE Find the intersection of constraints $\eta_{j'}$ and $\eta_{k'}$ 
and let its solution be $\bar{\x} = (\bar{x}_1,\bar{x}_2)$
\ENDIF
\IF{$k < m+2$}
\STATE $k \leftarrow k+1$
\ENDIF
\IF{$A'_{\eta_k,1}\bar{x}_1 + A'_{\eta_k,2}\bar{x}_2 > b'_{\eta_k}$}
\STATE $k' \leftarrow k$
\STATE Find the intersection of constraints $\eta_{j'}$ and 
$\eta_{k'}$ and let its solution be $\bar{\x} = (\bar{x}_1,\bar{x}_2)$
\ENDIF
\ENDWHILE
\STATE Report the optimal solution $\x' = \bar{\x}$ along 
with $\c'^{\Tr}\x'$, $\eta_{j'}$, and $\eta_{k'}$
\ENDIF
\end{algorithmic}
\label{SAlgorithm}
\end{algorithm}

Observe that finding an optimal basis to 2VLPs is critical when trying
to solve problem (\ref{twoDim}). This is because selecting two 
constraints that define an optimal solution but not an optimal basis 
(e.g. a 2VLP with an optimal degenerate solution) may result in a 
unchanged bases of the double pivot algorithm. This may result in 
unnecessary extra pivots and potentially, the algorithm may never 
terminate. For additional details on the benefits of the slope 
algorithm for degenerate linear programs, see 
\cite{ve18} and \cite{vitor18, vitor19, vitor22}. Furthermore, notice that degeneracy is not an issue for the slope algorithm when $\bar{x}_1=0$ 
or $\bar{x}_2=0$ since this problem has been resolved when 
selecting the two entering nonbasic variables. The following 
section presents the proposed double pivot method step by step and shows how the slope algorithm can be used to solve its subproblems.

\subsection{The Complete Double Pivot Algorithm}

The complete double pivot algorithm is provided in Algorithm \ref{mainAlg} and the following section presents the implementation details of the algorithm.

\begin{algorithm} {The Double Pivot Algorithm} \\
\begin{algorithmic}[1] 
\STATE Data: Matrix $\A$, vectors $\b$ and $\c$ 
\STATE Phase 1: Obtain an initial basic feasible solution $\x^0$ and its related partitions $\x_{B^0}$, $\x_{N^0}$, $\A_{B^0}$, $\A_{N^0}$, $\c_{B^0}$, and $\c_{N^0}$
\STATE Using (\ref{LUdecomp}), (\ref{LU1}), and (\ref{reducedCost}) calculate the reduced cost $\bar{\c}_{N^0}^{\Tr}=\c_{N^0}^{\Tr} - \c_{B^0}^{\Tr}  {\A_B^0} ^{-1} \A_N^0$ and determine the reduced cost vector index set $C^k$ using (\ref{RcostIndex})
\WHILE{$\min (\bar{\c}_{N^k}) < 0$}
\STATE Use Dantzig's pivoting rule to determine the entering variable $x_{\jmath_1^k}$, use (\ref{iniVar2}) to determine the
leaving variable $x_{\imath_1^k}$, then calculate the candidate 
$\x_B^{k+1}$ and the corresponding reduced cost $f_1$
\IF{the $|C^k| \ge 2$}
\STATE Use the longest step size rule for the set $C^k \setminus \{ \imath_1^k \}$ to determine the entering variable $x_{\jmath_2^k}$, use (\ref{iniVar3}) to determine the leaving variable $x_{\imath_2^k}$, then calculate the candidate 
$\x_B^{k+1}$ and the corresponding reduced cost  $f_2$
\STATE For the two entering variables $x_{\jmath_1^k}$ and $x_{\jmath_2^k}$, solve the two-dimensional LP problem (\ref{twoDim}) using Algorithm~\ref{SAlgorithm} to determine the leaving variable $x_{\imath_2^k}$ and $x_{\imath_2^k}$, then calculate the candidate $\x_B^{k+1}$ and the corresponding reduced cost $f_3$
\ENDIF
\IF{$|C^k| =1$}
\STATE Dantzig's pivoting pivot rule is used to update $B^k$ and $N^k$
\STATE LU decomposition is used for $\A_{B^k}$ and $\bar{\c}_{N^k}$
\ELSIF {$|C^k| \ge 2$}
\STATE The double pivoting rule is used to update $B^k$ and $N^k$
\STATE LU decomposition is used for $\A_{B^k}$ and $\bar{\c}_{N^k}$
\ENDIF
\STATE $k \leftarrow k+1$
\ENDWHILE
\end{algorithmic}
\label{mainAlg}
\end{algorithm}

\section{Implementation Details} \label{Imp_Det}

Some implementation details, which are important for improving the efficiency and robustness of Algorithm \ref{mainAlg}, are discussed in this section.

\subsection{Pre-Process and Post-Process}

Pre-process or pre-solver is a major factor that can significantly affect the numerical stability and computational efficiency of optimization algorithms. Many researchers have focused on this topic; for example, \cite{aa95}, \cite{bmw75}, \cite{mahajan10}, and \cite{yang17}. In this paper, the pre-process of \cite{yang17} is used. With the implementation described in the pre-process, the post-process is simple. The MATLAB code of the pre-process and post-process is available in the Netlib benchmark library (\textit{http://www.netlib.org/numeralgo/}) as part of  the \textit{na43} package \citep{yang17}.

\subsection{Update $\x_B^k$}

Once the updated base and $\A_{B^{k+1}}$ are available, one can update the basic feasible solution as $\x_{B^{k+1}}:=\A_{B^{k+1}}^{-1} \b$. However, for many Netlib benchmark problems, singular or nearly singular $\A_{B^{k+1}}$ are present. Thus, the computation of $\x_{B^{k+1}}$ using this formula may result in some negative components due to numerical errors. A better implementation is as follows. From (\ref{LP1}),

\begin{equation}
\x_{B^k}^{k+1}=\A_{B^k}^{-1}\b
-\A_{B^k}^{-1}\A_{N^k}\x_{N^k} 
= \bar{\b} - \bar{\A}_{(\jmath_1,\jmath_2)}
\x_{(\jmath_1,\jmath_2)} \ge \0.
\label{updatedXB1}
\end{equation}

\noindent If there is only one negative element in $\bar{\c}_{N^k}$, then

\begin{equation}
\x_{B^k}^{k+1}=\A_{B^k}^{-1}\b
-\A_{B^k}^{-1}\A_{N^k}\x_{N^k} 
= \bar{\b} - \bar{\A}_{\jmath_1}
x_{\jmath_1} \ge \0.
\label{updatedXB2}
\end{equation}

Observe that $\bar{\b}$, $\bar{\A}_{\jmath_1}$, and $\bar{\A}_{(\jmath_1,\jmath_2)}$ are computed by using LU decomposition for efficiency and robustness, and the details are discussed in Section \ref{avoid}. To obtain $\x_{B^{k+1}}^{k+1}$ from $\x_{B^k}^{k+1}$, one needs to remove the component of the leaving variable from $\x_{B^k}^{k+1}$ and insert the entering
variable into $\x_{B^k}^{k+1}$.

\subsection{Degenerate Solutions}

Even though the double pivoting rule solves all cycling problems of small size listed in Section \ref{cyclingTest}, there are some Netlib benchmark problems which stays in degenerate solutions for a long time. This is a sign that the double pivoting rule may still have difficult to handle all cycling problems. During the process of the code development, if a degenerate basic feasible solution moves to another degenerate basic feasible solution and the objective function is not improved after an iteration, then Bland's rule is applied \citep{bland77}. Still,
this implementation (using Bland's rule) has difficult to solve the Netlib 
benchmark problem QAP8 due to numerical errors. Inspired by Wolfe's perturbation 
method \citep{wolfe63}, the following perturbation is implemented: for a zero component in $\x_{B^k}$, a small positive perturbation to replace 
the zero component is introduced. This strategy clearly improves the stability of the code even when Bland's rule is not applied.

\subsection{LU Decomposition}

Normally for the simplex method, one would update the LU factorization with a rank-1 update (e.g., Forrest-Tomlin update \citep{ft72}), which saves substantial time. But MATLAB does not have such an implementation yet. For the purpose of comparing the computational efficiency of Dantzig's implementation and the double pivot algorithm, the general LU decomposition is implemented. This is not an optimal implementation, but it has the same effect on both algorithms; therefore, the efficiency comparison is still reasonable. However, the Forrest-Tomlin update will be considered in a future implementation.

The LU factorization can also be numerically unstable. This paper's implementation uses a dynamical permutation for improving the stability. Besides, when the diagonal elements of $\U$ become extremely small (close to zero), a positive perturbation is added to these elements. Computational experiments show the satisfactory numerical stability of the implementation.

\section{Numerical Tests} \label{Num_Tests}

Computational experiments were performed on an Intel$^{\circledR}$ Xeon$^{\circledR}$ E5-2670 2.60GHz 2 CPU/16 cores per node processor with 62.5GB of RAM per node. The version of MATLAB used was R2020a. The code that implements Algorithm \ref{mainAlg} has been extensively tested for many problems in 
some extreme cases, randomly generated problems, and benchmark problems. This section summarizes the test results. These test results have demonstrated
the computational merit of the double pivot algorithm.

\subsection{Test on the Klee-Minty Cube Problems} \label{klee}

The Klee-Minty cube and its variants have been used to show a serious drawback of popular simplex algorithms, such as steepest-edge pivoting rule, Dantzig's pivoting rule, most improvement pivoting rule, etc. i.e., all these algorithms need an exponential number of iterations, in the worst case, to find an optimal solution for specially designed Klee-Minty cube problems. This section provides the test results of the double pivot algorithm for three variants 
of the Klee-Minty cube \citep{greenberg97, ibrahima13, km11}. 

The first variant of the Klee-Minty cube is given by \cite{greenberg97}:

\begin{eqnarray}
\begin{array}{ll}
\min & \ \ -\sum_{i=1}^m 2^{m-i} x_i \\
\mbox{subject to} & 
\ \ \left[ 
\begin{array}{cccccc}
1 & 0 & 0 & \ldots & 0 & 0 \\
2^2 & 1 & 0 & \ldots & 0 & 0 \\
2^3 & 2^2  & 1 & \ldots & 0 & 0 \\
\vdots &  \vdots & \vdots &  \ddots & 0 & 0 \\
2^{m-1} &  2^{m-2} & 2^{m-3} & \ldots  & 1 & 0 \\
2^m &  2^{m-1} & 2^{m-2} & \ldots &  2^2  & 1 \\
\end{array}
\right]
\left[ \begin{array}{c}
x_1 \\ x_2 \\ \vdots \\  \vdots \\ x_{m-1} \\ x_m
\end{array} \right] 
\le 
\left[ \begin{array}{c}
5 \\ 25 \\ \vdots \\  \vdots \\ 5^{m-1} \\ 5^m
\end{array} \right] 
\\
& \ \ x_i \ge 0 \hspace{0.1in} i=1, \ldots, m.
\end{array}
\label{1stProblem}
\end{eqnarray}

\noindent The optimizer is $[0, \ldots, 0,5^m ]$ with an optimal objective function value of $-5^m$. The second variant of the Klee-Minty cube is given by \cite{ibrahima13}:

\begin{eqnarray}
\begin{array}{ll}
\min & \ \ -\sum_{i=1}^m 10^{m-i} x_i \\
\mbox{subject to} & 
\ \ \sum_{j=1}^{i-1} 10^{i-j} x_j +x_i \le 100^{i-1}
\hspace{0.1in} i=1, \ldots, m,  \\
& \ \ x_i \ge 0 \hspace{0.1in} i=1, \ldots, m.
\end{array}
\label{2ndProblem}
\end{eqnarray}

\noindent The optimizer is $[0, \ldots, 0,10^{2(m-1)}]$ with an optimal objective function value of $-10^{2(m-1)}$. The third variant of the Klee-Minty cube is given by \cite{km11}:

\begin{eqnarray}
\begin{array}{ll}
\min & \ \ -\sum_{i=1}^m   x_i \\
\mbox{subject to} & 
\ \ x_1 \le 1, \\
& \ \ 2 \sum_{i=1}^{k-1} x_i +x_k \le 2^k-1
\hspace{0.1in} k=2, \ldots, m, 
\\
& \ \ x_i \ge 0 \hspace{0.1in} i=1, \ldots, m.
\end{array}
\label{3rdProblem}
\end{eqnarray}

\noindent The optimizer is $[0, \ldots, 0,2^m-1]$ with an optimal objective function value of $-(2^m-1)$.

It is known that Dantzig's pivoting rule needs $2^m-1$ iterations to find an optimizer for these problems. But the proposed double pivot algorithm finds an optimizer in just one iteration \citep{yang20}. Therefore, the double pivot algorithm is able to solve these problems with size as large as $m=200$, which is impossible for Dantzig's algorithm because it needs about $2^{200} \ge 10^{60}$ iterations!

\subsection{Test on Randomly Generated Problems} \label{random}

This computational study also tested and compared Algorithm \ref{mainAlg} and Dantzig's pivoting algorithm using randomly generated problems. Notice that Algorithm \ref{mainAlg} is an improved version of the algorithm from \cite{yang20} because the slope algorithm is used for the two-dimensional LP problem. Therefore, the numerical results reported in this paper are better than the ones reported in \cite{yang20}. Computational experiments are carried out for randomly generated problems, which are obtained as follows: first, given the problem size $m$, a matrix $\M$ with uniformly distributed random entries between $[-0.5, 0.5]$ of dimension $m \times m$ and an identity matrix of dimension $m$ are generated. Therefore, $\A=[\M~~~ \I]$ is determined and the initial basic feasible solution is composed of the last $m$ columns. Then a positive vector $\b$ with uniformly distributed random entries between $[10, 11]$ of dimension $m$ and a vector $\c=(\c_1, \0)$ with $m$ $\c_1$'s entries uniformly distributed between $[-0.5, 0.5]$ are generated. 

Two MATLAB codes that implement Dantzig's pivoting algorithm and Algorithm \ref{mainAlg} are used to solve these randomly generated LP problems. Given the problem size $m$, $100$ random problems are generated and solved by using these two MATLAB codes. The average iteration number and average CPU time in seconds are 
obtained. The test results are presented in Table \ref{table1}. Notice that the numbers shown in Table \ref{table1} are the average of all 100 problems. Furthermore, CPU time is shown in seconds.

\begin{table}[h!]
\caption{Comparison Test for Dantzig's Pivoting Rule and the Double Pivot Algorithm for Randomly Generated Problems} \label{table1}
{\begin{tabular}{c|cc|cc|cc}
\hline
\textbf{Problem} & \multicolumn{2}{c|}{\textbf{Dantzig's Pivot}} & \multicolumn{2}{c|}{\textbf{Double Pivot}} & \multicolumn{2}{c}{\textbf{\% Improvement}}\\
\cline{2-3}
\cline{4-5}
\cline{6-7}
\textbf{Size $m$} & \textbf{Iter \#} & \textbf{CPU Time} & \textbf{Iter \#} & \textbf{CPU Time} & \textbf{Iter \#} & \textbf{CPU Time}\\
\hline
10&7&5.20E-4&5&1.59E-3&34\%&-205\%\\
100&165&8.62E-2&99&8.23E-2&40\%&4\%\\
1,000&17,334&4.58E+3&3,399&9.11E+2&80\%&80\%\\
10,000&160,741&6.97E+4&103,967&5.39E+4&35\%&23\%\\
\hline
\end{tabular}}
\end{table}

It is clear that the double pivot algorithm always uses few iterations to find an optimizer than Dantzig's algorithm does. However, for extremely small size problems, Dantzig's method uses less CPU time to find an optimizer than the double pivot algorithm does. As the problem size increases, the double pivot algorithm becomes more effective. For the problem size of $m=1000$, the double 
pivot algorithm improves performance by reducing about 80\% of the CPU 
time to find an optimizer. For the problem size of $m=10000$, the double pivot algorithm improves performance by reducing about 23\% of the CPU time to find an optimizer. For this paper, percentage improvement is computed as

\begin{equation}
\bigg(\displaystyle\frac{t_{\textsc{\tiny Dantzig}} - t_{\textsc{\tiny Double}}}{t_{\textsc{\tiny Dantzig}}}\bigg) \times 100\%
\end{equation}

\noindent where $t_{\textsc{\tiny Dantzig}}$ and $t_{\textsc{\tiny Double}}$
represents the number of pivots/CPU time required by Dantzig's implementation 
and the double pivot algorithm, respectively.

Observe that although most likely, $|C^k| \geq 2$, the cardinality of the randomly generated problems is much smaller than the ones from the Netlib benchmark problems. Therefore, the double pivot algorithm is more efficient than Dantzig's algorithm for the randomly generated problems with size from $100$ to
$10,000$. The double pivot algorithm is most efficient for problems with 
size of $1,000$. As the problem size increases further, $|C^k|$ becomes larger, and the double pivot algorithm becomes less efficient compared to Dantzig's algorithm.

\subsection{Test on Small Size Cycling Problems} \label{cyclingTest}

Notice that many Netlib benchmark problems have cycling issues. Before testing the large size benchmark Netlib problems, this computational study tested a collection of small size cycling problems listed in \cite{yang21}. This strategy helped the authors to develop a robust MATLAB code that appropriately implements the 
proposed algorithm to make it capable to avoid the cycling issue. As expected in Remark \ref{cycling}, the code for Algorithm \ref{mainAlg} performs well when degenerate solutions are encountered in these problems because the longest step size rule selects the entering variable with the largest value, which increases the chance to avoid the cycling issue. As a result, Algorithm \ref{mainAlg} solves all problems in this set!

\subsection{Test on Netlib Benchmark Problems}

The last set tested by this computational study is the Netlib benchmark library provided in \cite{bdgr95}. Both Dantzig's pivoting rule and the double pivot algorithm were tested for these problems. The size of these problems is
normally larger than the ones in the previous section, and finding the pivot that has the longest step size becomes extremely expensive. To reduce the computational cost, the pivot with longest step size  among candidates that satisfies the condition $\bar{c}_i <0.99\times\min \{\bar{\c} \}$ is selected. The test results are provided in Tables \ref{tableIteration} and \ref{tableIteration2}. Similarly, CPU time is shown in seconds. Moreover, due to the size of the table, results are split in two tables. It is clear that Dantzig's method uses less CPU time for almost all tested Netlib problems than the double pivot algorithm, even though the latter uses fewer iterations to find an optimizer. However, the double pivot algorithm still shows some merit on solving cycling problems. For problem {\tt DEGEN2}, even though the code for Dantzig's method implemented some safeguard tricks learned from the experience of the small size cycling problems, after more than one hundred million
iterations, Dantzig's pivoting rule still cannot find the solution, which is a sign that cycling has happened in this problem, but the double pivot algorithm finds the solution in only $2,350$ iterations.

\begin{landscape}
\begin{table}[h!]
\caption{Comparison Test for Dantzig's Pivoting Rule and the Double Pivot Algorithm for Netlib Benchmark Problems} \label{tableIteration}
{\begin{tabular}{c|cccc|cccc|cc}
\hline          
\multirow{2}{*}{\textbf{Problem}} & \multicolumn{4}{c|}{\textbf{Dantzig's Pivot}}  & \multicolumn{4}{c|}{\textbf{Double Pivot}} & \multicolumn{2}{c}{\textbf{\% Improvement}}\\ 
\cline{2-11}
& \textbf{Obj} & \textbf{Iter \#} & \textbf{CPU Time} & \textbf{Infea} & \textbf{Obj} & \textbf{Iter \#} & \textbf{CPU Time} & \textbf{Infea} & \textbf{Iter \#} & \textbf{CPU Time}\\
\hline
ADLITTLE&2.25495E+05&154&.09308&3.14578E-13&2.25495E+05&122&.26763&2.89899E-13&21\%&-188\%\\
AFIRO&-4.64753E+02&10&.04628&7.25485E-14&-4.64753E+02&8&.05018&3.82639E-14&20\%&-8\%\\
AGG&-3.59918E+07&469&.55888&1.44443E-10&-3.59918E+07&272&.84588&1.33442E-10&42\%&-51\%\\
AGG2&-2.02393E+07&560&.76756&1.28699E-10&-2.02393E+07&314&.88015&1.14947E-10&44\%&-15\%\\
AGG3&1.03121E+07&577&.79314&8.89876E-11&1.03121E+07&319&.91423&8.90887E-11&45\%&-15\%\\
BANDM&-3.01616E+02&581&.77758&2.11801E-13&-3.01616E+02&444&2.0567&4.01359E-13&24\%&-165\%\\
BEACONFD&3.20902E+04&79&.84003&2.38277E-13&3.20902E+04&44&.92386&1.39430E-13&44\%&-10\%\\
BLEND&-3.08121E+01&103&.07257&1.77005E-14&-3.08121E+01&97&.17077&1.49333E-14&6\%&-135\%\\
BNL1&1.97763E+03&3,577&9.21727&9.23917E-11&1.97763E+03&2,515&21.37082&9.15327E-11&30\%&-132\%\\
BNL2&1.77526E+03&13,220&84.74392&2.29648E-10&1.77526E+03&9,084&194.82839&2.58843E-10&31\%&-130\%\\
BRANDY&1.51851E+03&496&.28671&4.65831E-13&1.51851E+03&399&.95905&7.20008E-13&20\%&-235\%\\
DEGEN2&-&-&-&-&-1.43516E+03&2,350&22.04371&8.47990E-06&-&-\\
DEGEN3&-9.87287E+02&27,253&1236.393&1.71899E-09&-9.87290E+02&24,897&1753.88486&1.53891E-09&9\%&-42\%\\
FFFFF800&5.55680E+05&1,464&2.02042&6.69504E-10&5.55680E+05&1,101&6.59932&9.90309E-10&25\%&-227\%\\
ISRAEL&-8.96645E+05&875&.58372&2.52782E-10&-8.96645E+05&593&2.23322&3.73152E-10&32\%&-283\%\\
LOTFI&-2.52647E+01&286&.43982&7.76751E-10&-2.52647E+01&239&.66401&5.30878E-12&16\%&-51\%\\
MAROS\_R7&6.93586E+05&3,314&933.66064&3.16666E-10&6.93586E+05&1,991&821.81242&3.10953E-10&40\%&12\%\\
OSA\_07&5.35723E+05&6,378&303.95223&7.38573E-12&5.35723E+05&3,098&378.88184&5.69731E-12&51\%&-25\%\\
OSA\_14&1.10646E+06&13,111&1381.86663&3.22600E-12&1.10646E+06&6,787&2261.12124&1.62941E-10&48\%&-64\%\\
OSA\_30&2.14214E+06&27,623&4805.0524&3.39458E-10&2.14214E+06&14,618&8331.49749&3.63473E-10&47\%&-73\%\\
QAP8&2.03500E+02&44,232&2661.87395&3.07977E-09&2.03500E+02&45,125&3860.2692&3.47153E-09&-2\%&-45\%\\
SC50A&-6.45751E+01&22&.02811&3.57485E-14&-6.45751E+01&24&.04507&2.94039E-14&-9\%&-60\%\\
SC50B&-7.00000E+01&26&.02467&3.17764E-14&-7.00000E+01&22&.03733&8.32427E-14&15\%&-51\%\\
SC105&-5.22021E+01&52&.07303&1.03273E-13&-5.22021E+01&53&.12255&1.58714E-13&-2\%&-68\%\\
SC205&-5.22021E+01&91&.2418&2.36857E-12&-5.22021E+01&144&.49622&3.91708E-13&-58\%&-105\%\\
\hline
\end{tabular}}
\end{table}
\end{landscape}

\clearpage

\begin{landscape}
\begin{table}[h!]
\caption{Comparison Test for Dantzig's Pivoting Rule and the Double Pivot Algorithm for Netlib Benchmark Problems - Continued} \label{tableIteration2}
{\begin{tabular}{c|cccc|cccc|cc}
\hline          
\multirow{2}{*}{\textbf{Problem}} & \multicolumn{4}{c|}{\textbf{Dantzig's Pivot}}  & \multicolumn{4}{c|}{\textbf{Double Pivot}} & \multicolumn{2}{c}{\textbf{\% Improvement}}\\ 
\cline{2-11}
& \textbf{Obj} & \textbf{Iter \#} & \textbf{CPU Time} & \textbf{Infea} & \textbf{Obj} & \textbf{Iter \#} & \textbf{CPU Time} & \textbf{Infea} & \textbf{Iter \#} & \textbf{CPU Time}\\
\hline
SCAGR7&-2.08339E+06&201&.11251&2.83993E-12&-2.08339E+06&164&.33007&3.31881E-12&18\%&-193\%\\
SCAGR25&-1.45054E+07&1,043&2.39327&1.58450E-11&-1.45054E+07&1,016&7.77863&1.55783E-11&3\%&-225\%\\
SCFXM1&1.85315E+04&598&1.33481&1.04481E-09&1.85315E+04&462&1.96483&1.04478E-09&23\%&-47\%\\
SCFXM2&3.68898E+04&1,341&5.16887&3.65015E-12&3.68898E+04&958&7.57217&6.05719E-12&29\%&-46\%\\
SCFXM3&5.52456E+04&2,047&11.03826&5.74796E-12&5.52456E+04&1,547&16.92685&4.73435E-12&24\%&-53\%\\
SCRS8&9.04297E+02&616&3.07467&1.20544E-13&9.04297E+02&518&4.17995&2.89076E-14&16\%&-36\%\\
SCSD1&8.66667E+00&217&.05741&1.59068E-11&8.66667E+00&218&.37656&1.46231E-11&0\%&-556\%\\
SCSD6&5.05000E+01&528&.2355&3.13451E-11&5.05000E+01&401&1.05024&1.74808E-11&24\%&-346\%\\
SCSD8&9.05000E+02&1,242&1.43172&1.67105E-11&9.05000E+02&965&6.6592&1.48737E-11&22\%&-365\%\\
SCTAP1&1.41225E+03&691&.71354&3.09856E-10&1.41225E+03&534&2.06893&3.01522E-10&23\%&-190\%\\
SCTAP2&1.72481E+03&2,366&9.34237&5.88075E-14&1.72481E+03&2,228&33.37118&3.43825E-14&6\%&-257\%\\
SCTAP3&1.42400E+03&2,994&16.78474&6.60119E-10&1.42400E+03&2,537&46.74904&3.59637E-14&15\%&-179\%\\
SHARE1B&-7.01632E+04&513&.2243&1.78766E-10&-7.01632E+04&330&.52793&1.96484E-10&36\%&-135\%\\
SHARE2B&-3.58732E+02&166&.07706&5.41047E-11&-3.58732E+02&118&.16453&5.41023E-11&29\%&-113\%\\
SHIP04L&1.79315E+06&555&2.5496&3.34155E-14&1.79315E+06&258&12.89924&1.96751E-14&54\%&-406\%\\
SHIP04S&1.78245E+06&399&1.3644&1.37957E-14&1.78245E+06&190&5.48947&3.41211E-14&52\%&-302\%\\
SHIP08L&1.90734E+06&970&6.59985&2.41521E-14&1.90734E+06&650&25.11337&2.47146E-14&33\%&-281\%\\
SHIP08S&1.88314E+06&541&2.554&2.29440E-14&1.88314E+06&324&5.16255&2.29741E-14&40\%&-102\%\\
SHIP12L&1.46555E+06&1,211&12.38697&3.89742E-14&1.46555E+06&791&28.34739&3.11866E-14&35\%&-129\%\\
SHIP12S&1.46060E+06&600&4.28729&3.34777E-14&1.46060E+06&352&10.30063&4.01002E-14&41\%&-140\%\\
STOCFOR1&-4.11320E+04&57&.12437&1.60489E-12&-4.11320E+04&39&.14181&1.60489E-12&32\%&-14\%\\
STOCFOR2&-3.90244E+04&1,733&39.92682&1.46442E-11&-3.90244E+04&1,303&51.91037&1.56429E-11&25\%&-30\%\\
STOCFOR3&-3.99768E+04&14,453&2571.60678&4.58717E-11&-3.99768E+04&14,590&4985.93169&4.44197E-11&-1\%&-94\%\\
TRUSS&4.58816E+05&12,656&99.91333&2.38932E-11&4.58816E+05&7,315&181.77245&2.37071E-11&42\%&-82\%\\
\hline
\end{tabular}}
\end{table}
\end{landscape}

\section{Conclusions} \label{Conc}

This paper proposed a double pivot algorithm that combines the pivoting rule of \cite{yang20} and the slope algorithm of \cite{ve18} and \cite{vitor18, vitor19, vitor22}. This paper's research implemented the algorithm as a MATLAB code. Computational experiments tested the MATLAB code using a small set of cycling problems, three variants of the Klee-Minty problems, randomly generated LP problems, and Netlib benchmark test problems. A computational study compared the performances of the double pivot algorithm and the code that implemented Dantzig's pivoting algorithm. The test results show that the double pivot algorithm performs better than Dantzig's pivoting algorithm for large size randomly generated problems, while the latter performs better for small size randomly generated problems and
Netlib benchmark test problems.

\section*{Acknowledgement}
A portion of this work was partially funded by the National Science Foundation 
(NSF) under the EPSCoR research program -- Grant N$^{\circ}$ OIA-1557417. 
This work was completed utilizing the Holland Computing Center of the University
of Nebraska, which receives support from the Nebraska Research Initiative.

\end{document}